\documentclass[12pt]{article}
\usepackage{hyperref,amsfonts,amssymb}
\usepackage{amsfonts,amssymb}
\usepackage{mathtools}
\usepackage{rotating}

\def\C{\mathbb{C}}

\def\Z{\mathbb{Z}}

\def\bq{ \begin{equation} }
\def\eq{ \end{equation} }
\def\ben{ \begin{eqnarray} }
\def\en{ \end{eqnarray} }

\def\frac#1#2{{#1\over #2}}

\def\on#1#2{\mathop{\vbox{\ialign{##\crcr\noalign{\kern2pt}
$\scriptstyle{#2}$\crcr\noalign{\kern2pt\nointerlineskip}
\kern-2pt$\hfil\displaystyle{#1}\hfil$\crcr}}}\limits}

\begin{document}

\title{Functional equations in algebra}
\author{Boris Feigin,  Alexander Odesskii}
   \date{}
\vspace{-20mm}
   \maketitle
\vspace{-7mm}

\medskip

\begin{abstract}

\medskip

We study flat deformations of quotients of a polynomial algebra in a class of graded commutative associative algebras. Functional equations and their solutions in terms of theta functions play important role in these studies. An analog of this theory in a fermionic case is also briefly discussed.  

\medskip

\end{abstract}

\newpage
\tableofcontents
\newpage

\section{Introduction}
 
An ideal $I\subset\C[x_1,...,x_N]$ is called monomial if it is generated by a set of monomials in $\C[x_1,...,x_N]$. Monomial ideals often appear in the crossroad of commutative algebra, combinatorics and algebraic geometry \cite{cox}. In this paper we study flat deformations\footnote{See Definition 4.1.} of quotients by monomial ideals. Some of these deformations might play important role in representation theory of infinite dimensional Lie algebras \cite{kac} and in the theory of vertex algebras \cite{fr}.

\subsection{Some experimental data} 

Let $P^{(N)}=\C[x_{-N},...,x_N]$ where $N\gg0$. Let $I_0\subset P^{(N)}$ be a monomial ideal generated by 
$$\{x_ix_{i+1},~-N\leq i\leq N-1\}.$$
Let us study homogeneous flat deformations in this case. The quotient $P^{(N)}/I_0$ admits a monomial basis 
\begin{equation}\label{mb}
x_{i_1}^{r_1}...x_{i_m}^{r_m},~~~i_1+1<i_2,~i_2+1<i_3,...,i_{m-1}+1<i_m 
\end{equation}
where $-N\leq i_1,...,i_m\leq N,~0<r_1,...,r_m$.

Let $I\subset P^{(N)}$ be an ideal generated by 
\begin{equation}\label{i0}
\{x_ix_{i+1}+a_1x_{i-1}x_{i+2}+a_2x_{i-2}x_{i+3}+...,~-N\leq i\leq N-1\} 
\end{equation}
where $a_1,a_2,...\in\C$ are parameters. We consider $I$ as a deformation of $I_0$. This deformation is said to be flat if the algebra $P^{(N)}/I$ has the same basis (\ref{mb}) as the algebra 
$P^{(N)}/I_0$. We want to find all parameters $a_1,a_2,...$ such that this deformation is flat. To obtain constraints for these parameters we proceed as follows. Consider a monomial\footnote{We assume that $N$ is much larger than all indexes appear in this computation.}
$$x_ix_{i+1}x_{i+2}\in P^{(N)}/I.$$
We can transform it in two different ways:
$$x_ix_{i+1}x_{i+2}=-(a_1x_{i-1}x_{i+2}+a_2x_{i-2}x_{i+3}+...)x_{i+2}=-a_1x_{i-1}x_{i+2}^2-a_2x_{i-2}x_{i+2}x_{i+3}-...$$
or
$$x_ix_{i+1}x_{i+2}=-x_i(a_1x_ix_{i+3}+a_2x_{i-1}x_{i+4}+...)=-a_1x_i^2x_{i+3}-a_2x_{i-1}x_ix_{i+4}-...$$
Let us replace all products of the form $x_jx_{j+1}$ by $-a_1x_{j-1}x_{j+2}-a_2x_{j-2}x_{j+3}-...$ everywhere in the r.h.s. of both expressions for $x_ix_{i+1}x_{i+2}$. For example, 
we replace $x_{i+2}x_{i+3}$ in the first expression and $x_{i-1}x_i$ in the second expression. Continuing in this way we obtain two expressions for $x_ix_{i+1}x_{i+2}$ in terms of 
monomials of the form (\ref{mb}). We assume that these monomials (\ref{mb}) are linear independent. Equating coefficients of each monomial of the two expressions thus obtained for $x_ix_{i+1}x_{i+2}$ we get constraints for $a_1,a_2,...$. General solution of these constraints is as follows. Let $I_q,~q\in\C$ be a family of ideals generated by 
\begin{equation}\label{reli}
\sum_{\alpha\in\Z}(-1)^{\alpha}q^{\frac{\alpha(3\alpha-1)}{2}}x_{i+3\alpha}x_{i+1-3\alpha},~i\in\Z. 
\end{equation}
If a deformation $P^{(N)}/I$ is flat, then $I=I_q$ for some $q\in\C$.
On the other hand, it is proved in this paper (Theorem 5.1) that indeed, the algebra $P^{(N)}/I$ is a flat deformation of $P^{(N)}/I_0$.

Similar computations can be done for other algebras with quadratic monomial relations, see Conjectures 5.1, 5.2 for description of results of these experiments.

\subsection{Functional equations}

Let us introduce a generating function for coefficients $a_1,a_2,...$ in (\ref{i0})
$$f(z_1,z_2)=z_1+z_2+\sum_{j=1}^{\infty} a_j(z_1^{-j}z_2^{j+1}+z_1^{j+1}z_2^{-j}).$$
We have $f(z_2,z_1)=f(z_1,z_2)$ and $f(\lambda z_1,\lambda z_2)=\lambda f(z_1,z_2).$ It turns out that a deformation $P^{(N)}/I$ is flat iff the function $f(z_1,z_2)$ satisfies a functional equation
$$z_1^2g\Big(\frac{z_2z_3}{z_1^2}\Big)f(z_2,z_3)+z_2^2g\Big(\frac{z_3z_1}{z_2^2}\Big)f(z_3,z_1)+z_3^2g\Big(\frac{z_1z_2}{z_3^2}\Big)f(z_1,z_2)=0$$
where 
$$g(z)=\sum_{j\in\Z}b_jz^j$$
is another function. This and similar functional equations play important role in our study of flat deformations, see Section 5 for details.

\subsection{Connections with Representation Theory and Vertex Algebras}

It is convenient to set $N=\infty$ in the previous considerations, this can be done formally in terms of projective limit in a category of graded algebras, see Section 2 for details.

Introduce a generating function
\begin{equation}\label{gf}
G(z)=\sum_{i=-\infty}^{\infty}z^ix_i 
\end{equation}
where $x_i,~i\in\Z$ are generators and $z$ is a formal parameter. It is known \cite{f2} that the commutative associative algebra with generators $x_i,~i\in\Z$ and relations $G(z)^{k+1}=0$ 
can be used for a realization of integrable representations of affine Lie algebra $\hat{sl}_2$ of level $k$. Moreover, the algebra with relations 
\begin{equation}\label{Gk}
G(z)G(tz)...G(t^kz)=0
\end{equation}
plays a similar role for a quantum analog \cite{dr} of $\hat{sl}_2$. 

Consider for example the case 
\begin{equation}\label{G2}
G(z)G(tz)=0. 
\end{equation}
Substituting (\ref{gf}) into (\ref{G2}), expanding and equating to zero coefficients of all powers of $z$ we obtain the relations
$$x_i^2+(t+t^{-1})x_{i-1}x_{i+1}+(t^2+t^{-2})x_{i-2}x_{i+2}+...=0,$$
$$(1+t)x_ix_{i+1}+(t^{-1}+t^2)x_{i-1}x_{i+2}+(t^{-2}+t^3)x_{i-2}x_{i+3}+...=0.$$
On the other hand, we can deform relations (\ref{G2}) in a similar way as we have deformed the relations $x_ix_{i+1}=0$ into  relations (\ref{reli}). We obtain the following relations
\begin{equation}\label{G2d}
\sum_{\alpha\in\Z}(-1)^{\alpha}q^{\frac{\alpha(3\alpha-1)}{2}}G(t^{3\alpha}z)G(t^{1-3\alpha}z)=0. 
\end{equation}
Substituting the expression (\ref{gf}) into (\ref{G2d}) and equating to zero coefficients of all powers of $z$ we obtain a family of algebras. We believe that this is a flat deformation of an algebra with relations $x_i^2=0,~x_ix_{i+1}=0$, see Conjecture 5.1 for details. 

In Theorem 5.2 we have constructed a flat deformation of the algebra with relations 
$$x_ix_{i+1}...x_{i+k}=0.$$ 
Probably a similar deformation of relations (\ref{Gk}) gives an interesting family 
of algebras. We believe that the study of these and similar algebras, particularly, the study of flat deformations of commutative algebras constructed in \cite{f1} may have applications in representation theory and the theory of vertex algebras.

\subsection{Plan of the paper and main results}

In Section 2 we explain our formalism of commutative associative algebras possessing infinite expressions. We also introduce our notations for theta functions.

In Section 3 we explain so-called functional realization of polynomials and other algebras in terms of Laurent series. This is our main tool in obtaining flatness constraints in terms of  
functional equations.

In Section 4 we obtain necessary and sufficient conditions for flatness of homogeneous deformations of an algebra with relations $x_ix_{i+1}=0$ (Theorem 4.1) and, more generally, with relations 
$x_ix_{i+1}...x_{i+k-1}=0$ with fixed $k>1$ (Theorem 4.2).

In Section 5 we construct flat deformations of an algebra with relations $x_ix_{i+1}=0$ (Theorem 5.1) and, more generally, with relations 
$x_ix_{i+1}...x_{i+k-1}=0$ with fixed $k>1$ (Theorem 5.2). We  construct a conjectural flat deformation of the algebra with relations $x_i^2=0,~x_ix_{i+1}=0$ (Conjecture 5.1). 
We also briefly discuss a fermionic case which is similar to bosonic one. We suggest a conjectural flat deformation of an algebra with relations $\xi_i\xi_{i+1}=0$ (Conjecture 5.2) and, more generally, of an algebra with relations $\xi_i\xi_{i+1}...\xi_{i+k-1}=0$ (Conjecture 5.3).

In Section 6 we discuss possible applications of our results and directions of future research.

\section{Notations and conventions}

{\bf ~~~Algebra $A$.~} Consider a polynomial algebra $P=\C[...,x_{-2},x_{-1},x_0,x_1,...]$ generated by infinitely many variables $x_i,~i\in\Z$. This is a $\Z^2$-graded associative commutative algebra: we have
$P=\oplus_{n,l}P_{n,l}$ where $P_{n,l}$ is spanned by monomials $x_{i_1}^{a_1}...x_{i_m}^{a_m}$ with $a_1i_1+...+a_mi_m=n,~a_1+...+a_m=l$. We want to extend this algebra in the following way: 
let $A_{n,l}$ consist of infinite formal sums of monomials in $P_{n,l}$ and let $A=\oplus_{n,l}A_{n,l}$. For example, $\sum_{i\in\Z}c_ix_{3+i}x_{-i}\in A_{3,2}$ where $c_i\in\C$ are 
arbitrary coefficients. It is clear that multiplication of such infinite sums is well defined and so $A$ is an associative commutative algebra. Another way of defining $A$ is as follows. Let $P^{(N)}=\C[x_{-N},...,x_N]$. Define homomorphisms $\phi_N:~P^{(N+1)}\to P^{(N)}$ by 
$$\phi_N(x_{-N-1})=\phi_N(x_{N+1})=0~ \text{and}~ \phi_N(x_i)=x_i~ \text{if} ~-N\leq i\leq N.$$
The algebra $A$ is the projective limit in the category of $\Z^2$-graded algebras corresponding to this chain of homomorphisms.

{\bf Homogeneous ideals of the algebra $A$.~} In this paper we study quotients of the $\Z^2$-graded algebra $A$ by homogeneous ideals. When we define an ideal in $A$ we use the following convention. Let $g_i,~i\in\Z$ be a set of homogeneous elements of $A$. By an ideal generated by $\{g_i,~i\in\Z\}$ we always mean the set of infinite formal sums $\sum_{i}p_ig_i\in A$ where 
$p_i\in A$. Notice that each of these formal sums should belongs to a direct sum of finitely many homogeneous components $A_{n,l}$. Quotients of $A$ by such ideals can also be defined as projective limits. Let us illustrate this by the following

{\bf Example 2.1.} Let $Q_N$ be a $\Z^2$-graded algebra with generators $x_i,~-N\leq i\leq N$ and defining relations 
$$\sum_{\substack{\alpha+\beta=i,\\-N\leq \alpha,\beta\leq N}}x_{\alpha}x_{\beta}=0,~~~i\in\Z.$$
We assume that $\text{deg}(x_i)=(i,1)\in\Z^2$. For each $N=1,2,...$ there exists a homomorphism $\phi_N:~Q_{N+1}\to Q_N$ such that 
$$\phi(x_{-N-1})=\phi(x_{N+1})=0~ \text{and}~ \phi(x_i)=x_i~ \text{if} ~-N\leq i\leq N.$$
Let $Q$ be the projective limit in the category of $\Z^2$-graded algebras corresponding to this chain of homomorphisms.

{\em In the sequel we will say in this and similar examples that the algebra $Q$ is the quotient of the algebra $A$ by an ideal generated by 
$$\sum_{\substack{\alpha+\beta=i,\\\alpha,\beta\in\Z}}x_{\alpha}x_{\beta}=0,~~~i\in\Z.$$}

{\bf Dimensions of graded components and cut-off ideals $C_N$.~} Let $I(q)$ be a family of graded ideals of the algebra $A$ where $q\in\C$ is a parameter. In this paper we study how dimensions 
of graded components of a quotient $A/I(q)$ depend on $q$. For example, we would like to say that for certain family of ideals $\dim (A/I(q))_{m,l}$ does not depend on $q$ or that $\dim (A/I(q))_{m,l}\leq \dim (A/I(0))_{m,l}$ etc. A technical problem here is that  $\dim A_{m,l}=\infty$ and also  $\dim (A/I(q))_{m,l}=\infty$ for all families of ideals we are going to study. To overcome this problem we introduce cut-off ideals $C_N$ of the algebra $A$. By definition $C_N$ is the ideal generated by $x_{N+1},x_{N+2},..$. In other words, we cut-off all $x_i$ with $i>N$. 
It is clear that $\dim (A/C_N)_{m,l}<\infty$ for all $N$. Let us introduce a notation $\dim_N$. By definition
$$\dim_N Q_{m,l}=\dim (Q/C_N)_{m,l}$$
where $Q$ is a quotient of $A$. In the sequel when we use a notation $\dim_N$ we always mean that $N\gg 0$.

{\bf Theta functions.~} Fix a complex number $\tau$ such that $\text{Im}\tau>0$. A theta function in one variable of degree one is defined by
$$\theta(u)=\sum_{k\in\Z}(-1)^ke^{2\pi i(ku+\frac{k(k-1)}{2}\tau)}.$$
In this paper we use multiplicative form of this and other theta functions. Introduce new variables $z=e^{2\pi iu},~q=e^{2\pi i\tau}$ and denote $\theta(u)$ in these variables by $g(z)$. We have
$$g(z)=\sum_{k\in\Z}(-1)^kz^iq^{\frac{k(k-1)}{2}}.$$
It is known \cite{th} that the function $g(z)$ can be written as an infinite product
$$g(z)=(1-z)\prod_{i=1}^{\infty}(1-q^i)(1-q^iz)(1-q^iz^{-1}).$$
See \cite{th} for this and other facts about theta functions.

\section{Functional realization}

Recall that $A=\oplus_{n,l}A_{n,l}$ is a $\Z^2$-graded  algebra, we assume $\text{deg}(x_i)=(i,1)$.  Sometimes it will be convenient to identify $\oplus_nA_{n,l}$ with a vector space of symmetric Laurent formal series in $l$ variables. Let us define isomorphisms $\psi_l:~\oplus_kA_{n,l}\to\C((z_1^{\pm1},...,z_l^{\pm1}))^{S_l}$ by 
\begin{equation}\label{fr}
\psi_l(x_{i_1}...x_{i_l})=\sum_{\sigma\in S_l}z_{\sigma(1)}^{i_1}...z_{\sigma(l)}^{i_l} 
\end{equation}
where $S_l$ stands for symmetric group. For example $\psi_1(x_i)=z_1^i$ and $\psi_2(x_ix_j)=z_1^iz_2^j+z_1^jz_2^i$. The product of Laurent series $A_{\alpha}\otimes A_{\beta}\to A_{\alpha+\beta}$ can be written under such identification as follows
\begin{equation}\label{fr1}
fg(z_1,...,z_{\alpha+\beta})=\frac{1}{\alpha!\beta!}\sum_{\sigma\in S_{\alpha+\beta}}f(z_{\sigma(1)},...,z_{\sigma(\alpha)})g(z_{\sigma(\alpha+1)},...,z_{\sigma(\alpha+\beta)}) 
\end{equation}
where $f(z_1,...,z_{\alpha})\in\C((z_1^{\pm1},...,z_{\alpha}^{\pm1}))^{S_{\alpha}}$ and $g(z_1,...,z_{\beta})\in\C((z_1^{\pm1},...,z_{\beta}^{\pm1}))^{S_{\beta}}$. 

For example, if $f(z_1),~g(z_1)\in\C((z_1^{\pm1}))$, then 
$$fg(z_1,z_2)=f(z_1)g(z_2)+f(z_2)g(z_1)$$
and if $f(z_1)\in\C((z_1^{\pm1})),~g(z_1,z_2)\in\C((z_1^{\pm1},z_2^{\pm1}))^{S_2}$, then $$fg(z_1,z_2,z_3)=f(z_1)g(z_2,z_3)+f(z_2)g(z_1,z_3)+f(z_3)g(z_1,z_2).$$

We will refer to the identification (\ref{fr}) and to the product formula (\ref{fr1}) as to a functional realization of $A$. 

\section{Flat deformations of algebras with monomial relations}

Let $I_0\subset A$ be an ideal in $A$ generated by a set of monomials $m_{\alpha}$. The quotient $A/I_0$ is a $\Z^2$-graded algebra because each monomial $m_{\alpha}$ is homogeneous. Let $q\in\C$ be a parameter, and let $m_{\alpha}(q)$ be a family of homogeneous elements of $A$ such that $m_{\alpha}(0)=m_{\alpha}$. We say that the ideal $I_q$ is a deformation of $I_0$ and the quotient  $A/I_q$ is a deformation of $A/I_0$.

{\bf Lemma 4.1.} If $q\in\C$ is in generic position, then $\dim_N (A/I_q)_{n,l}\leq \dim_N (A/I_0)_{n,l}$ for all $n,l$. More precisely, for each $n,l$ there exists $\varepsilon_{n,l}>0$ such that $\dim_N (A/I_q)_{n,l}\leq \dim_N (A/I_0)_{n,l}$ for $0<|q|<\varepsilon_{n,l}$. 

{\bf Proof.} Follows from the usual semi-continuity of rank arguments. Indeed, $\dim_N (I_q)_{n,l}$ is equal to the rank of a certain matrix depending on $q$. Therefore, it does not depend on $q$ if $q$ is in generic position and can be only smaller for special values of $q$. $\Box$

{\bf Definition 4.1.} A deformation $A/I_q$ is flat if dimensions of graded components $(A/I_q)_{n,l}$ and $(A/I_0)_{n,l}$ are equal for all $n,l$. 

In this Section we discuss criteria of flatness of such deformations for various set of monomials $m_{\alpha}$. We will give necessary and sufficient conditions of flatness in the case 
$m_{\alpha}=x_{\alpha}x_{\alpha+1}...x_{\alpha+k},~\alpha\in\Z$ where $k>0$ is fixed. Specific examples of deformations will be discussed in the next Section.

\subsection{An algebra with quadratic relations}

Let $I_0$ be the ideal of the algebra $A$ generated by monomials $y_i=x_ix_{i+1},~i\in\Z$. In this Subsection we discuss homogeneous flat deformations of the quotient $A/I_0$. It is clear that the following relations hold 
\begin{equation}\label{rel}
x_iy_{i+1}-x_{i+2}y_i=0,~~~i\in\Z. 
\end{equation}
We are going to prove that a deformation of the algebra $A/I_0$ is flat iff the relations (\ref{rel}) can also be deformed.

Let $B$ be the quotient of the algebra $A$ by an ideal generated by 
\begin{equation}\label{z}
z_i=x_ix_{i+1}+a_1x_{i-1}x_{i+2}+a_2x_{i-2}x_{i+3}+... 
\end{equation}
 where $a_1,a_2,...\in\C$ are fixed constants. After a change of variable $x_j\to q^{j^2}x_j$ where $q\ne0$ 
is a parameter, we obtain an isomorphic algebra $B(q)$, the quotient of the algebra $A$ by the ideal generated by 
\begin{equation}\label{zq}
z_i(q)=x_ix_{i+1}+q^4a_1x_{i-1}x_{i+2}+q^{12}a_2x_{i-2}x_{i+3}+... 
\end{equation}
It is clear that $B(q)\cong B(1)=B$ if $q\ne 0$. On the other hand, $B(0)$ is a quotient of $A$ by the ideal generated by the monomials $z_i(0)=y_i=x_ix_{i+1}$ and is not isomorphic to $B$. The algebra $B(q)$ is a deformation of $B(0)$, the question is for which $a_1,a_2,...$ it is flat.  

{\bf Theorem 4.1.} $B(q)$ is a flat deformation of $B(0)$ iff $z_i(q)$ satisfy
\begin{equation}\label{reldef}
x_iz_{i+1}(q)-x_{i+2}z_i(q)=\sum_{\beta\in\Z} g_{\beta}(q)x_{j-2\beta}z_{j+1+\beta}(q),~~~i\in\Z. 
\end{equation}
where $g_{\beta}(q)$ are some holomorphic functions  and $g_{\beta}(0)=0$. 

Proof is based on combinatorial properties of the ideal $I_0$ generated by $y_i,~i\in\Z$ in a polynomial algebra.

Let $M$ be the quotient of the polynomial algebra $P=\C[x_i,~i\in\Z]$ by the ideal generated by $y_i=x_ix_{i+1},~i\in \Z$. It is clear that $M$ becomes a $\Z^2$-graded algebra if we assume  $\text{deg}(x_i)=(i,1)\in\Z^2$. The algebra $M$ has a basis $x_{i_1}^{a_1}x_{i_2}^{a_2}...x_{i_m}^{a_m}$ where $i_1,...,i_m\in\Z$, $0< a_1,...,a_m$ and $i_1+1<i_2,i_2+1<i_3,...,i_{m-1}+1<i_m$. 

Introduce an auxiliary polynomial 
algebra $\bar{P}=\C[x_i,\bar{y}_i,~i\in\Z]$. Define a homomorphism $\phi:~\bar{P}\to P$ by 
$$\phi(x_i)=x_i,~\phi(\bar{y}_i)=x_ix_{i+1}.$$

{\bf Definition 4.2.} A monomial $x_{i_1}^{a_1}...x_{i_m}^{a_m}\bar{y}_{j_1}^{b_1}...\bar{y}_{j_n}^{b_n}$ with $i_1<...<i_m,~j_1<...<j_n\in\Z$, $0<a_1,...,a_m,b_1,...,b_n$ is reduced if $i_1+1< i_2,...,i_{m-1}+1< i_m$ and $j_{\alpha}\ne i_{\beta}+1$ for $\alpha=1,...,n,~\beta=1,...,m$.

{\bf Remark 4.1.} In other words, a monomial is reduced if it is not divisible by $x_ix_{i+1}$ and $x_i\bar{y}_{i+1}$ for all $i\in\Z$.

{\bf Definition 4.3.} Let $u,v$ be monomials in the algebra $\bar{P}$. A monomial $v$ is a reduction of a monomial $u$ if $v$ can be obtained from $u$ by a sequence of the following transformations

{\bf 1.} Replacing $x_ix_{i+1}$ by $\bar{y}_i$ for some $i\in\Z$. In particular, $\bar{y}_i$ is a reduction of $x_ix_{i+1}$.

{\bf 2.} Replacing $x_i\bar{y}_{i+1}$ by $x_{i+2}\bar{y}_i$ for some $i\in\Z$. In particular, $x_{i+2}\bar{y}_i$ is a reduction of $x_i\bar{y}_{i+1}$.

{\bf Remark 4.2.} It is clear that a monomial $u$ is reduced iff no transformations 1 or 2 can be applied to $u$.

{\bf Lemma 4.2.} Let $u,v$ be monomials in the algebra $\bar{P}$. If a monomial $v$ is a reduction of a monomial $u$, then $\phi(u)=\phi(v)$. In particular, for any monomial $u\in \bar{P}$ the set of its reductions is finite.

{\bf Proof.} It is clear that $\phi(x_ix_{i+1})=\phi(\bar{y}_i)=x_ix_{i+1}$ and $\phi(x_i\bar{y}_{i+1})=\phi(x_{i+2}\bar{y}_i)=x_ix_{i+1}x_{i+2}$, so each of the transformations 1, 2 from Definition 4.3 preserves the image 
of a monomial under the map $\phi$. Therefore, any sequence of these transformations also preserves the image. 

Notice that for any monomial  $x_{l_1}...x_{l_t}\in P$ the set of monomials $v\in \bar{P}$ such that $\phi(v)=x_{l_1}...x_{l_t}$ is finite. Therefore, the set of reductions of any monomial $u\in\bar{P}$ is also finite.                                $\Box$

{\bf Lemma 4.3.} For any monomial $u$ in the algebra $\bar{P}$ there exists a unique reduced monomial $v$ such that $v$ is a reduction of $u$. 

{\bf Proof.} For an arbitrary monomial $v=x_{i_1}^{a_1}...x_{i_m}^{a_m}~\bar{y}_{j_1}^{b_1}...\bar{y}_{j_n}^{b_n}\in\bar{P}$ we denote $\text{deg}_1(v)=a_1+...+a_m$ and $\text{deg}_2(v)=b_1j_1+...+b_nj_n$. Let 
$$\text{Red}(u)=\{v\in\bar{P};~v~\text{is a monomial and is a reduction of}~u\}.$$
The set $\text{Red}(u)$ is finite by Lemma 4.2. Note that\footnote{$u$ is obtained from $u$ by an empty sequence of transformations from Definition 4.3.} 
$u\in\text{Red}(u)$ so the set $\text{Red}(u)$ is not empty. Take a monomial $v\in\text{Red}(u)$ and apply transformations 1 and 2 to it until it is possible. Notice that each transformation 1 decreases $\text{deg}_1(v)$ by one and each transformation 2 preserves $\text{deg}_1(v)$ and  decreases $\text{deg}_2(v)$ by 1. Since $\text{Red}(u)$ is finite, after finitely many applications of transformations 1 and 2 we come to an element $w\in\text{Red}(u)$ such that neither transformation 1 nor 2 can be applied to $w$. This means that $w$ is not divisible by $x_ix_{i+1}$ and $x_i\bar{y}_{i+1}$ for all $i\in\Z$. Therefore, $w$ is reduced by Remark 4.1 and we have proved existence.

To prove uniqueness, suppose that $v,v^{\prime}$ are both reductions of $u$ and are both reduced. It follows from Lemma 4.2 that $\phi(v)=\phi(v^{\prime})=\phi(u)$ and the statement $v=v^{\prime}$ 
follows from Lemma 4.4 below. $\Box$

{\bf Lemma 4.4.} If monomials $v,v^{\prime}\in\bar{P}$ are both reduced and $\phi(v)=\phi(v^{\prime})$, then $v=v^{\prime}$.

{\bf Proof.} Let $v=x_{i_1}^{a_1}...x_{i_m}^{a_m}~\bar{y}_{j_1}^{b_1}...\bar{y}_{j_n}^{b_n},~v^{\prime}=x_{i^{\prime}_1}^{a^{\prime}_1}...x_{i^{\prime}_{m^{\prime}}}^{a^{\prime}_{m^{\prime}}}~\bar{y}_{j^{\prime}_1}^{b^{\prime}_1}...\bar{y}_{j^{\prime}_{n^{\prime}}}^{b^{\prime}_{n^{\prime}}}$ where $i_1<...<i_m,~j_1<...<j_n\in\Z$, $0<a_1,...,a_m,b_1,...,b_n$ and $i^{\prime}_1<...<i^{\prime}_{m^{\prime}},~j^{\prime}_1<...<j^{\prime}_{n^{\prime}}\in\Z$, $0<a^{\prime}_1,...,a^{\prime}_{m^{\prime}},b^{\prime}_1,...,b^{\prime}_{n^{\prime}}$. Assume that $v,v^{\prime}$ are both reduced and $\phi(v)=\phi(v^{\prime})$. Assume that $v\ne v^{\prime}$ and we have chosen a pair $v,v^{\prime}$ with the smallest possible $g=a_1+...+a_m+b_1+...+b_n$. We have $i_1\ne i^{\prime}_1$ and $j_1\ne j^{\prime}_1$, otherwise we can divide both $v,v^{\prime}$ by $x_{i_1}$ or $y_{j_1}$ and decrease $g$. Let $i_1<i^{\prime}_1$ (otherwise just interchange $v$ and $v^{\prime}$). In this case we have $j^{\prime}_1=i_1$ because $\phi(v)=\phi(v^{\prime})$ is divisible by $x_{i_1}$. On the other hand, 
$\phi(\bar{y}_{j^{\prime}_1})=\phi(\bar{y}_{i_1})=x_{i_1}x_{i_1+1}$ and therefore $\phi(v)$ is divisible by $x_{i_1+1}$. This is not possible however because $v$ is reduced by assumption and so 
$i_2\ne i_1+1$ and $j_a\ne i_1+1$ for all $a=1,...,n$. $\Box$

{\bf Remark 4.3.} The statements of the above Lemmas can be reformulated less formally as follows: all relations between $x_i,y_i\in P$ are consequences of the relations (\ref{rel}).

Let us prove Theorem 4.1. If a deformation $B(q)$ is flat, then in particular $\dim_N B(q)_{m,3}=\dim_N B(0)_{m,3}$ for all $m$. This means that the relations (\ref{rel}) can be deformed and we have (\ref{reldef}) for some functions $g_{\beta}(q)$. Indeed, if the relations (\ref{rel}) can not be deformed, then the cubic part of the ideal $I_q$ becomes larger than the cubic part of the ideal $I_0$.

On the other hand, suppose that (\ref{reldef}) holds. It follows from Lemma 4.1 that $\dim_N B(q)_{n,l}\leq \dim_N B(0)_{n,l}$ for $q\ne 0$ because all algebras $B(q)$ with $q\ne 0$ are isomorphic.
 Therefore, the dimensions of graded components of the ideal $I_q$ generated by $z_i(q)$ for $q\ne 0$ can be either the same or larger than these at $q=0$. However, it follows from Lemma 4.3 that        the ideal generated by $z_i(q)$ for  $q\ne 0$ can not become larger than the one for $q=0$ because any element of this ideal of the form $x_{i_1}^{a_1}...x_{i_m}^{a_m}z_{j_1}(q)^{b_1}...z_{j_n}(q)^{b_n}$ with $i_1<...<i_m,~j_1<...<j_n\in\Z$, $0<a_1,...,a_m,b_1,...,b_n$ can be written as a linear combination of reduced monomials of such form. To show this let us write (\ref{zq})  in the form
\begin{equation}\label{tr1}
x_ix_{i+1}=z_i(q)+O(q) 
\end{equation}
Let us also write (\ref{reldef}) in the form
\begin{equation}\label{tr2}
x_iz(q)_{i+1}=x_{i+2}z(q)_i+O(q) 
\end{equation}
Given a monomial $u=x_{i_1}^{a_1}...x_{i_m}^{a_m}z_{j_1}(q)^{b_1}...z_{j_n}(q)^{b_n}\in I(q)$ we can make the following transformations without changing $u\in A$

 ${\bf1^{\prime}.~}$ Replacing $x_ix_{i+1}$ by $z_i(q)+O(q)$ for some $i\in\Z$.

 ${\bf2^{\prime}.~}$ Replacing $x_iz_{i+1}$ by $x_{i+2}z_i(q)+O(q)$ for some $i\in\Z$.
 
 These transformations coincide modulo $q$ with transformations 1 and 2 from Definition 4.3  if we replace $z(q)_i$ by $y_i$. Lemma 4.3 shows that $u$ is equal to a reduced monomial modulo $q$. Write 
 $u=u^{(r)}+q\sum_{\alpha}s_{\alpha} v_{\alpha} + O(q^2)$ where $u^{(r)}$ is a reduced monomial, $s_{\alpha}\in\C$ and $v_{\alpha}$ are some monomials. After that we apply the same procedure to each monomial $v_{\alpha}$ writing it as $v^{(r)}_{\alpha}+O(q)$ where $v^{(r)}_{\alpha}$ are reduced. Continuing in this way we write $u$ as a linear combination of reduced monomials.
 
 Since each monomial from the ideal $I_q$ can be written as a linear combination of reduced monomials, then $\dim_N (I_q)_{n,l}$ can not be greater than the number of reduced monomials of degree $(n,l)$. But this number of reduced monomials is equal to $\dim_N (I_0)_{n,l}$ by Lemma 4.3.      $\Box$

\subsection{Algebras with higher degree relations}

Fix an integer $k>1$. Let $I_0$ be the ideal of the algebra $A$ generated by monomials $y_i=x_ix_{i+1}...x_{i+k},~i\in\Z$. In this Subsection\footnote{In the case $k=1$ we have quadratic relations. This case was considered in the previous Subsection.} we discuss homogeneous flat deformations of the quotient $A/I_0$. It is clear that the following relations hold 
\begin{equation}\label{relk}
x_iy_{i+1}-x_{i+k+1}y_i=0,~~~i\in\Z. 
\end{equation}
We are going to prove that a deformation of the algebra $A/I_0$ is flat iff the relations (\ref{relk}) can also be deformed.

Let $B$ be the quotient of the algebra $A$ by an ideal generated by 
\begin{equation}\label{zk}
z_i=\sum_{\substack{\alpha_1,...,\alpha_{k+1}\in\Z,\\ \alpha_1\leq...\leq\alpha_{k+1},\\ \alpha_1+...+\alpha_{k+1}=\frac{k(k+1)}{2}}} a_{\alpha_1,...,\alpha_{k+1}}x_{i+\alpha_1}x_{i+\alpha_2}...x_{i+\alpha_{k+1}},~~~i\in\Z
\end{equation}
 for some fixed constants $a_{\alpha_1,...,\alpha_{k+1}}\in\C$ where without loss of generality we assume $a_{0,1,...,k}=1$. After change of variable $x_j\to q^{j^2}x_j$ where $q\ne0$ 
is a parameter, we obtain an isomorphic algebra $B(q)$, the quotient of the algebra $A$ by the ideal generated by 
\begin{equation}\label{zqk}
z_i(q)=\sum_{\substack{\alpha_1,...,\alpha_{k+1}\in\Z,\\ \alpha_1\leq...\leq\alpha_{k+1},\\ \alpha_1+...+\alpha_{k+1}=\frac{k(k+1)}{2}}} q^{\alpha_1^2+...+\alpha_{k+1}^2-1^2-...k^2}a_{\alpha_1,...,\alpha_{k+1}}x_{i+\alpha_1}x_{i+\alpha_2}...x_{i+\alpha_{k+1}},~~~i\in\Z
\end{equation}
It is clear that $B(q)\cong B(1)=B$ if $q\ne 0$. Assume that the r.h.s. of (\ref{zqk}) can be written as $x_ix_{i+1}...x_{i+k}+O(q)$. In this case $B(0)$ is a quotient of $A$ by the ideal generated by the monomials $z_i(0)=y_i=x_ix_{i+1}...x_{i+k}$ and is not isomorphic to $B$. On the other hand, $B(q)$ is a deformation of $B(0)$, the question is for which $a_{\alpha_1,...,\alpha_{k+1}}$ it is flat.  

{\bf Theorem 4.2.} $B(q)$ is a flat deformation of $B(0)$ iff $z_i(q)$ satisfy
\begin{equation}\label{reldefk}
x_iz_{i+1}(q)-x_{i+k+1}z_i(q)=\sum_{\beta\in\Z} g_{\beta}(q)x_{j-(k+1)\beta}z_{j+\beta}(q),~~~i\in\Z. 
\end{equation}
where $g_{\beta}(q)$ are some holomorphic functions  and $g_{\beta}(0)=0$. 

Proof is similar to the proof of Theorem 4.1. We outline the main steps leaving details to the reader.

Let $M$ be the quotient of the polynomial algebra $P=\C[x_i,~i\in\Z]$ by the ideal generated by $y_i=x_ix_{i+1}...x_{i+k},~i\in \Z$. Introduce an auxiliary polynomial 
algebra $\bar{P}=\C[x_i,\bar{y}_i,~i\in\Z]$. Define a homomorphism $\phi:~\bar{P}\to P$ by 
$$\phi(x_i)=x_i,~\phi(\bar{y}_i)=x_ix_{i+1}...x_{i+k}.$$

{\bf Definition 4.4.} A monomial $x_{i_1}^{a_1}...x_{i_m}^{a_m}\bar{y}_{j_1}^{b_1}...\bar{y}_{j_n}^{b_n}$   is reduced if it is not divisible by $x_ix_{i+1}...x_{i+k}$ and $x_i\bar{y}_{i+1}$ for all $i\in\Z$.

{\bf Definition 4.5.} Let $u,v$ be monomials in the algebra $\bar{P}$. A monomial $v$ is a reduction of a monomial $u$ if $v$ can be obtained from $u$ by a sequence of the following transformations

{\bf 1.} Replacing $x_ix_{i+1}...x_{i+k}$ by $\bar{y}_i$ for some $i\in\Z$. In particular, $\bar{y}_i$ is a reduction of $x_ix_{i+1}...x_{i+k}$.

{\bf 2.} Replacing $x_i\bar{y}_{i+1}$ by $x_{i+k+1}\bar{y}_i$ for some $i\in\Z$. In particular, $x_{i+k+1}\bar{y}_i$ is a reduction of $x_i\bar{y}_{i+1}$.

With these Definitions Lemmas 4.2, 4.3, 4.4 still hold, proofs are similar. Proof of Theorem 4.2 is based on these Lemmas and is similar to proof of Theorem 4.1.  $\Box$

\section{Flat deformations and functional equations}

\subsection{Algebras with quadratic relations}

\bigskip

{\bf 1. From flatness conditions to functional equations.}

Let $I_0$ be an ideal in $A$ generated by monomials 
$$y^0_{1,i}=x_ix_{i+k_1},...,y^0_{l,i}=x_ix_{i+k_l},~i\in\Z$$ 
where $0\leq k_1<...<k_l$ are fixed. There exist relations between $x_i$ and 
$y^0_{a,j}$ of the form 
\begin{equation}\label{relq}
\sum_{\substack{j\in\Z\\1\leq a\leq l}}p_{a,j,s}x_{i+s-2j-k_a}y^0_{a,i+j}=0,~~~i\in\Z 
\end{equation}
where $p_{a,j,s}\in\{0,1,-1\}$. For example, we have relations $x_iy^0_{a,i+k_b}-x_{i+k_a+k_b}y^0_{b,i}=0$. Let $D_s$ be dimension of the vector space of relations of the form (\ref{relq}) with 
given $s$. Note that  $D_s>0$ for some $s$ unless $l=1,~k_1=0$. 

{\bf Example 5.1.} Let $I_0$ be generated by $y^0_i=x_ix_{i+1}$. We have relations $x_iy^0_{i+1}-x_{i+2}y^0_i=0$ so $D_3=1$.

{\bf Example 5.2.} Let $I_0$ be generated by $y^0_{1,i}=x_i^2,~y^0_{2,i}=x_ix_{i+1}$. We have relations 
$$x_iy^0_{2,i}-x_{i+1}y^0_{1,i}=0,~~~x_iy^0_{1,i+1}-x_{i+1}y^0_{2,i}=0,~~~x_iy^0_{2,i+1}-x_{i+2}y^0_{2,i}=0$$
so $D_1=D_2=D_3=1$.

Let $I$ be an ideal in $A$ generated by 
\begin{equation}\label{def}
y_{a,i}=x_ix_{i+k_a}+\sum_{j\in\Z\setminus \{0,-k_a\}}u_{a,j}x_{i-j}x_{i+k_a+j}
\end{equation}
where $u_{a,j}\in\C$. The question is for which parameters $u_{a,j}$ the quotient algebra $A/I$ has the same size as $A/I_0$. More precisely, we want to have $\dim_N(A/I)_{m,l}=\dim_N(A/I_0)_{m,l}$ for all $m,l$ and $N\gg 0$. In this case we say that $A/I$ is a flat deformation\footnote{See Section 4 for precise definition.} of $A/I_0$. In particular we want $\dim_N(A/I)_{m,3}=\dim_N(A/I_0)_{m,3}$ which means that relations 
\begin{equation}\label{rel1}
\sum_{\substack{j\in\Z\\1\leq a\leq l}}v_{a,j,s}x_{i+s-2j-k_a}y_{a,i+j}=0,~~~i\in\Z 
\end{equation}
should hold where $v_{a,j,s}\in\C$. Moreover, dimension of vector space of such relations should be equal to $D_s$. 
Substituting the expressions (\ref{def}) into (\ref{rel1}) and equating to zero all coefficients of monomials in $x_k$ we obtain constraints for $u_{a,j},v_{a,j,s}\in\C$ which should be solved in order to classify all flat deformations. It is much easier however to do this in terms of generating functions for the coefficients $u_{a,j},v_{a,j,s}$. Let us do this using functional realization described in Section 3.

Introduce generating functions for the sequences $u_{a,j}$, $v_{a,j,s}$
\begin{equation}\label{f1}
f_{a}(z_1,z_2)=z_1^{k_a}+z_2^{k_a}+\sum_{j\in\Z\setminus \{0,-k_a\}} u_{a,j}(z_1^{-j}z_2^{k_a+j}+z_1^{k_a+j}z_2^{-j}),~~~g_{a,s}(z)=\sum_{j\in\Z}v_{a,j,s}z^j.
\end{equation}
It is clear that the following equations hold
\begin{equation}\label{hom1}
f_a(\lambda z_1,\lambda z_2)=\lambda^{k_a}f_a(z_1,z_2). 
\end{equation}
and
\begin{equation}\label{fa2}
f_a(z_1,z_2)=f_a(z_2,z_1). 
\end{equation}

{\bf Lemma 5.1.} Relations (\ref{rel1}) are equivalent to the following functional equation
\begin{equation}\label{fa1}
\sum_{a=1}^l (z_1^{s-k_a}g_{a,s}\Big(\frac{z_2z_3}{z_1^2}\Big)f_a(z_2,z_3)+z_2^{s-k_a}g_{a,s}\Big(\frac{z_3z_1}{z_2^2}\Big)f_a(z_3,z_1)+z_3^{s-k_a}g_{a,s}\Big(\frac{z_1z_2}{z_3^2}\Big)f_a(z_1,z_2))=0.
\end{equation}
Moreover, if functions $f_a(z_1,z_2)$ are fixed, then for each $s$ dimension of vector space of functions $g_{a,s}(z)$ satisfying (\ref{fa1}) is equal to $D_s$.

{\bf Proof.} Let\footnote{See Section 3 for definition of $\psi_2$.} 
$$f_{a,i}(z_1,z_2)=\psi_2(y_{a,i})=z_1^iz_2^{i+k_a}+z_1^{i+k_a}z_2^i+\sum_{j\in\Z\setminus \{0,-k_a\}} u_{a,j}(z_1^{i-j}z_2^{i+k_a+j}+z_1^{i+k_a+j}z_2^{i-j}).$$
It is clear that  we can write 
\begin{equation}\label{fi}
f_{a,i}(z_1,z_2)=z_1^iz_2^if_a(z_1,z_2)
\end{equation}
Computing the functional realization of the l.h.s. of (\ref{rel1}) using formula (\ref{fr1}) for multiplication we obtain
$$\sum_{\substack{j\in\Z\\1\leq a\leq l}}v_{a,j,s}(z_1^{i+s-2j-k_a}f_{a,i+j}(z_2,z_3)+z_2^{i+s-2j-k_a}f_{a,i+j}(z_3,z_1)+z_3^{i+s-2j-k_a}f_{a,i+j}(z_1,z_2))=0.$$
Using (\ref{fi}) we rewrite this as 
\begin{equation}\label{fa}
\sum_{\substack{j\in\Z\\1\leq a\leq l}}v_{a,j,s}(z_1^{s-2j-k_a}z_2^jz_3^jf_a(z_2,z_3)+z_2^{s-2j-k_a}z_3^jz_1^jf_a(z_3,z_1)+z_3^{s-2j-k_a}z_1^jz_2^jf_a(z_1,z_2))=0. 
\end{equation}
Using generating functions $g_{a,s}(z)$ for the sequences $v_{a,j,s}$ we can rewrite (\ref{fa}) as (\ref{fa1}). $\Box$

\bigskip

{\bf 2. Flat deformations of the algebra with relations $x_ix_{i+1}=0$.}

Let $I_0$ be generated by $y^0_i=x_ix_{i+1}$. In this case $D_3=1$, see Example 1. Equation (\ref{fa1}) takes the form
\begin{equation}\label{fa11}
z_1^2g_{1,0}\Big(\frac{z_2z_3}{z_1^2}\Big)f_1(z_2,z_3)+z_2^2g_{1,0}\Big(\frac{z_3z_1}{z_2^2}\Big)f_1(z_3,z_1)+z_3^2g_{1,0}\Big(\frac{z_1z_2}{z_3^2}\Big)f_1(z_1,z_2)=0 
\end{equation}

{\bf Lemma 5.2.} The functions $g_{1,0}(z)=g(z),~f_1(z_1,z_2)$ defined by 
\begin{equation}\label{fg}
\begin{array}{ccc}
g(z)=(1-z)\prod_{i=1}^{\infty}(1-q^i)(1-q^iz)(1-q^iz^{-1})=\sum_{i\in\Z}(-1)^iz^iq^{\frac{i(i-1)}{2}}, \\
\\
f_1(z_1,z_2)=\frac{z_2g\big(\frac{z_1^2}{z_2^2}\big)}{g\big(\frac{z_1}{z_2}\big)}\prod_{i=1}^{\infty}(1-q^i)=\sum_{i\in\Z}(-1)^iq^{\frac{i(3i-1)}{2}}(z_1^{3i}z_2^{1-3i}+z_1^{1-3i}z_2^{3i})
\end{array}
\end{equation}
satisfy equation (\ref{fa11}). Here $q\in\C$ and $|q|<1$.

{\bf Proof.} It is clear that the Laurent series in (\ref{fg}) is convergent if $|q|<1$.  Substituting the first expression for $f_1(z_1,z_2)$ from (\ref{fg}) into (\ref{fa11}) we obtain
\begin{equation}\label{fag}
t_1^2g\Big(\frac{t_2}{t_1^2}\Big)\frac{g(t_2^2)}{g(t_2)}+t_1t_2^2g\Big(\frac{t_1}{t_2^2}\Big)\frac{g\big(\frac{1}{t_1^2}\big)}{g\big(\frac{1}{t_1}\big)}+t_2g(t_1t_2)\frac{g\big(\frac{t_1^2}{t_2^2}\big)}{g\big(\frac{t_1}{t_2}\big)}=0
\end{equation}
where $t_1=\frac{z_1}{z_3},~t_2=\frac{z_2}{z_3}$.

One can verify by direct computations the following identities
\begin{equation}\label{id}
 g(z^{-1})=-z^{-1}g(z),~~~g(1)=0,~~~g(qz)=-z^{-1}g(z).
\end{equation}
Let $R(t_1,t_2)$ be the l.h.s. of (\ref{fag}). The function $R(t_1,t_2)$ possesses a Laurent series expansion in $t_1,t_2$. One can verify using (\ref{id}) that 
$$R(qt_1,t_2)=q^{-1}t_1^{-4}t_2^2R(t_1,t_2),~~~R(t_1,qt_2)=q^{-1}t_1^2t_2^{-4}R(t_1,t_2)$$
and
$$R(t_2^{-1},t_2)=R(t_2^2,t_2)=R(\pm t_2^{1/2},t_2)=R(-1,t_2)=R(-t_2,t_2)=0.$$
There is no nonzero Laurent series in $t_1,~t_2$ satisfying these properties. $\Box$

{\bf Theorem 5.1.} Let $B_q$ be an algebra generated by $x_i,~i\in\Z$ with defining relations 
\begin{equation}\label{relB}
\sum_{\alpha\in\Z}(-1)^{\alpha}q^{\frac{\alpha(3\alpha-1)}{2}}x_{i+3\alpha}x_{i+1-3\alpha}=0,~i\in\Z. 
\end{equation}
Then $B_q$ is a flat deformation of $B_0$.

{\bf Proof.} Follows from Lemmas 5.1, 5.2 and Theorem 4.1. Notice that if we define 
$$y_i=\sum_{\alpha\in\Z}(-1)^{\alpha}q^{\frac{\alpha(3\alpha-1)}{2}}x_{i+3\alpha}x_{i+1-3\alpha}$$
then in the algebra $A$ we have
$$\sum_{\beta\in\Z}(-1)^{\beta}q^{\frac{\beta(\beta-1)}{2}}x_{j-2\beta+2}y_{j+\beta}=0,~j\in\Z.~~~\Box$$

\bigskip

{\bf 3. Deformations of the algebra with relations $x_i^2=x_ix_{i+1}=0$.} 

Let $I_0$ be generated by $y^0_{1,i}=x_i^2,~y^0_{2,i}=x_ix_{i+1}$. We have relations 
$$x_iy^0_{2,i}-x_{i+1}y^0_{1,i}=0,~~~x_iy^0_{1,i+1}-x_{i+1}y^0_{2,i}=0,~~~x_iy^0_{2,i+1}-x_{i+2}y^0_{2,i}=0$$
so $D_1=D_2=D_3=1$.

Define a function $f(t)$ by
$$f(t)=\sum_{i\in\Z}(-1)^i\Big(t^{6i+1}+\frac{1}{t^{6i+1}}\Big)q^{\frac{i(3i+1)}{2}}.$$
where $|q|<1$.

{\bf Conjecture 5.1.} Let $I$ be an ideal in the algebra $A$ generated by 
$$y_{1,i}=\sum_{k\in\Z}f(t^{2k})\tilde{q}^{k^2}x_{i-k}x_{i+k},~~~y_{2,i}=\sum_{k\in\Z}f(t^{2k+1})\tilde{q}^{k^2+k}x_{i-k}x_{i+k+1},~~~i\in\Z$$
where $t,q,\tilde{q}\in\C$, $t\ne0$. The algebra $A/I$ is a flat deformation of the algebra $A/I_0$.

{\bf Remark 5.1.} To support this conjecture we have verified that the corresponding functional equations (\ref{fa1}) hold for certain functions $g_{a,s}(z)$.

{\bf Remark 5.2.} The class of isomorphism of the algebra $A/I$ does not depend on $\tilde{q}$ if $\tilde{q}\ne0$.

{\bf Remark 5.3.} Let $\tilde{q}=1$. Introduce a generating function $G(z)=\sum_{i\in\Z}z^ix_i$, a formal Laurent series in a parameter $z$. Relations in the algebra $A/I$ can be written as
$$\sum_{k\in\Z}(-1)^kq^{\frac{k(3k+1)}{2}}G\Big(\frac{z}{t^{6k+1}}\Big)G(zt^{6k+1})=0.$$
In particular, if $q=0$, then we have $G(\frac{z}{t})G(zt)=0$.

\subsection{Algebras with higher degree relations}

Let $B_{k,0}$, $2\leq k$ be the quotient of the polynomial algebra $A=\C[x_i,~i\in\Z]$ by the ideal generated by $y^0_i=x_ix_{i+1}...x_{i+k-1}\in A$. It is clear that $B_{k,0}$ becomes a $\Z^2$-graded algebra if we assume  $\text{deg}(x_i)=(i,1)\in\Z^2$. 

We want to study flat deformations of $B_{k,0}$ in the class of $\Z^2$-graded algebras of the type described in Section 4.2. Recall that the following relations hold
\begin{equation}\label{relh}
x_iy^0_{i+1}-x_{i+k}y^0_i=0,~~~i\in\Z.  
\end{equation}
Suppose that a deformed $\Z^2$-graded algebra $B$ has relations 
\begin{equation}\label{defh}
y_i=\sum_{\substack{\alpha_1,...,\alpha_k\in\Z,\\ \alpha_1\leq...\leq\alpha_k,\\ \alpha_1+...+\alpha_k=\frac{k(k-1)}{2}}} a_{\alpha_1,...,\alpha_k}x_{i+\alpha_1}x_{i+\alpha_2}...x_{i+\alpha_k},~~~i\in\Z
\end{equation}
for some $a_{\alpha_1,...,\alpha_k}\in\C$ where without loss of generality we assume $a_{0,1,...,k-1}=1$. In order to obtain a flat deformation we need also to deform the relations between relations (\ref{relh}). Suppose that after deformation we have
\begin{equation}\label{relh1}
\sum_{\beta\in\Z} b_{\beta}x_{j+k\beta}y_{j+1-\beta}=0,~~~j\in\Z. 
\end{equation}
Substituting the expressions (\ref{defh}) into (\ref{relh1}) and equating to zero all coefficients of monomials in $x_k$ we obtain constraints for $a_{\alpha_1,...,\alpha_k},b_{\beta}\in\C$ which should be solved in order to classify all flat deformations. It is much easier however to do this in term of generating functions for the coefficients $a_{\alpha_1,...,\alpha_k},b_{\beta}\in\C$. Let us do this using functional realization described in Section 3. Introduce generating functions for the sequences $a_{\alpha_1,...,\alpha_k},b_{\beta}$
\begin{equation}\label{fh1}
f(z_1,...,z_k)=\sum_{\substack{\alpha_1,...,\alpha_k\in\Z,\\ \alpha_1\leq...\leq\alpha_k,\\ \alpha_1+...+\alpha_k=\frac{k(k-1)}{2}\\\sigma\in S_k}} a_{\alpha_1,...,\alpha_k}z_{\sigma(1)}^{\alpha_1+1}z_{\sigma(2)}^{\alpha_2+1}...z_{\sigma(k)}^{\alpha_k+1} 
\end{equation}
and
$$g(z)=\sum_{\beta\in\Z}b_{\beta}z^{\beta}.$$
Note that $f(z_1,...,z_k)$ is symmetric and satisfies the homogeneity condition
\begin{equation}\label{hom}
f(\lambda z_1,...,\lambda z_k)=\lambda^{\frac{k(k+1)}{2}}f(z_1,...,z_k). 
\end{equation}

{\bf Lemma 5.3.} Relations (\ref{relh1}) are equivalent to the following functional equation
\begin{equation}\label{fah1}
g\Big(\frac{z_1^k}{z_2...z_{k+1}}\Big)f(z_2,...,z_{k+1})+...+g\Big(\frac{z_{k+1}^k}{z_1...z_k}\Big)f(z_1,...,z_k)=0.
\end{equation}

{\bf Proof.} Let\footnote{See Section 3 for definition of $\psi_k$.} 
$$f_i(z_1,...,z_k)=\psi_k(y_i)=\sum_{\substack{\alpha_1,...,\alpha_k\in\Z,\\ \alpha_1\leq...\leq\alpha_k,\\ \alpha_1+...+\alpha_k=\frac{k(k-1)}{2}\\\sigma\in S_k}} a_{\alpha_1,...,\alpha_k}z_{\sigma(1)}^{i+\alpha_1}z_{\sigma(2)}^{i+\alpha_2}...z_{\sigma(k)}^{i+\alpha_k}.$$
It is clear that  we can write 
\begin{equation}\label{fih}
f_i(z_1,...,z_k)=z_1^{i-1}...z_k^{i-1}f(z_1,...,z_k) 
\end{equation}
Computing the functional realization of the l.h.s. of (\ref{relh1}) using the formula (\ref{fr1}) for multiplication we obtain
$$\sum_{\beta\in\Z} b_{\beta}(z_1^{j+k\beta}f_{j+1-\beta}(z_2,...,z_{k+1})+...+z_{k+1}^{j+k\beta}f_{j+1-\beta}(z_1,...,z_k))=0,~~~j\in\Z.$$
Using (\ref{fih}) we rewrite this as 
\begin{equation}\label{fah}
\sum_{\beta\in\Z} b_{\beta}(z_1^{k\beta}z_2^{-\beta}...z_{k+1}^{-\beta}f(z_2,...,z_{k+1})+...+z_1^{-\beta}...z_k^{-\beta}z_{k+1}^{k\beta}f(z_1,...,z_k))=0.
\end{equation}
Using generating function $g(z)$ for the sequence $b_{\beta}$ 
we can rewrite (\ref{fah}) as (\ref{fah1}) 
where the function $f(z_1,...,z_k)$ is symmetric and satisfies the equation (\ref{hom}). $\Box$

Define functions $g(z)$, $f_{n,k}(z_1,...,z_n)$ by 
\begin{equation}\label{fun}
\begin{array}{ccc}
g(z)=(1-z)\prod_{i=1}^{\infty}(1-q^i)(1-q^iz)(1-q^iz^{-1})=\sum_{i\in\Z}(-1)^iz^iq^{\frac{i(i-1)}{2}}, \\
\\
f_{n,k}(z_1,...,z_n)=(-1)^n\sum_{\sigma\in S_n} \frac{g(z_{\sigma(1)}^k)g(z_{\sigma(1)}^{-1}z_{\sigma(2)}^k)g(z_{\sigma(1)}^{-1}z_{\sigma(2)}^{-1}z_{\sigma(3)}^k)...g(z_{\sigma(1)}^{-1}...z_{\sigma(n-1)}^{-1}z_{\sigma(n)}^k)}{g(z_{\sigma(1)}^{-1})g(z_{\sigma(1)}^{-1}z_{\sigma(2)}^{-1})g(z_{\sigma(1)}^{-1}z_{\sigma(2)}^{-1}z_{\sigma(3)}^{-1})...g(z_{\sigma(1)}^{-1}...z_{\sigma(n)}^{-1})}
\end{array}
\end{equation}
where $|q|<1$ and $n,k=1,2,...$. For example,
$$f_{1,k}(z)=-\frac{g(z^k)}{g(z^{-1})}~~~\text{and}~~~f_{2,k}(z_1,z_2)=\frac{g(z_1^k)g(z_1^{-1}z_2^k)}{g(z_1^{-1})g(z_1^{-1}z_2^{-1})}+\frac{g(z_2^k)g(z_2^{-1}z_1^k)}{g(z_2^{-1})g(z_1^{-1}z_2^{-1})}.$$
The following properties of the functions $g,~f_{n,k}$ can be verified by direct calculation
\begin{equation}\label{gpr}
g(qz)=-z^{-1}g(z);~~~g(z)=0 ~\text{iff} ~z=q^i,~i\in\Z, 
\end{equation}
\begin{equation}\label{fpr}
f_{n,k}(qz_1,z_2,...,z_n)=(-1)^{k-1}q^{-\frac{(k+1)(k-2)}{2}}z_1^{-k^2+1}z_2^{k+1}...z_n^{k+1}f_{n,k}(z_1,...,z_n). 
\end{equation}
Furthermore, $f_{n,k}(z_1,...,z_n)$ is a symmetric function in $z_1,...,z_n$ and satisfies the following recurrence relations
\begin{equation}\label{req}
f_{n+1,k}(z_1,...,z_{n+1})=-\sum_{i=1}^{n+1}\frac{g(z_1^{-1}...z_i^k...z_{n+1}^{-1})}{g(z_1^{-1}...z_{n+1}^{-1})}f_{n,k}(z_1,...,\widehat{z_i},...,z_{n+1}). 
\end{equation}

It follows from (\ref{gpr}), (\ref{fun}) that 
\begin{equation}\label{rest1}
f_{1,k}(1)=k,~~~f_{1,k}(\varepsilon)=0,~~~f_{1,k}(\eta)=(-1)^{i-1}\eta^iq^{-\frac{i(i-1)}{2}} 
\end{equation}
where $\varepsilon^k=q^i,~\eta^{k+1}=q^i,~i\in\Z$ and $\varepsilon,\eta\ne 1$. More generally, we have
\begin{equation}\label{rest2}
f_{n+1,k}(z_1,...,z_n,1)=(k-n)f_{n,k}(z_1,...,z_n) 
\end{equation}
and
\begin{equation}\label{rest3}
f_{n+1,k}(z_1,...,z_n,\varepsilon)=-z_1^{(k+1)i}f_{n,k}(\varepsilon z_1,z_2,...,z_n)-...-z_n^{(k+1)i}f_{n,k}(z_1,...z_{n-1},\varepsilon z_n), 
\end{equation}
\begin{equation}\label{rest4}
\begin{array}{cc}
f_{n+1,k}(z_1,...,z_n,\eta)=(-1)^{i-1}z_1^i...z_n^i\eta^iq^{-\frac{i(i-1)}{2}} f_{n,k}(z_1,...,z_n)-\\ \\ -z_1^{ki}f_{n,k}(\eta z_1,z_2,...,z_n)-...-z_n^{ki}f_{n,k}(z_1,...z_{n-1},\eta z_n) 
\end{array}
\end{equation}
where $\varepsilon^k=q^i,~\eta^{k+1}=q^i,~i\in\Z$ and $\varepsilon,\eta\ne 1$. This follows either directly from (\ref{fun}) or using induction by $n$ from (\ref{req}).

{\bf Lemma 5.4.} The functions $f_{n,k}(z_1,...,z_n)$ are holomorphic if $z_1,...,z_n\in\C\setminus\{0\}$ and $|q|<1$. In particular, $f_{n,k}(z_1,...,z_n)$ can be represented as Laurent series in $z_1,...,z_n$.

{\bf Proof.} Let us use induction in $n$. The case $n=1$ is clear by (\ref{fun}). Suppose that $f_{n,k}(z_1,...,z_n)$ is holomorphic. We need to prove that the r.h.s. of (\ref{req}) is 
holomorphic as well. Let 
$$Res(z_1,...,z_n)=-\sum_{i=1}^{n+1}g(z_1^{-1}...z_i^k...z_{n+1}^{-1})f_{n,k}(z_1,...,\widehat{z_i},...,z_{n+1})\Big |_{z_{n+1}=z_1^{-1}...z_n^{-1}}.$$
We need to prove that $Res(z_1,...,z_n)=0$ identically. It follows from the induction hypothesis that $Res(z_1,...,z_n)$ is holomorphic. It is also clear that $Res(z_1,...,z_n)$ is symmetric and satisfies the identity
\begin{equation}\label{per}
Res(qz_1,...,z_n)=q^{-k(k+1)}z_1^{-2k^2-2k}z_2^{-k^2-k}...z_n^{-k^2-k}Res(z_1,...,z_n). 
\end{equation}
It follows from (\ref{rest1}), (\ref{rest2}), (\ref{rest3}) that $Res(z_1,...,z_n)$ is divisible by
$$\frac{g(z_1^k)...g(z_n^k)g(z_1^{-k}...z_n^{-k})~g(z_1^{k+1})...g(z_n^{k+1})g(z_1^{-k-1}...z_n^{-k-1})}{g(z_1)...g(z_n)g(z_1^{-1}...z_n^{-1})}.$$
Let
$$H(z_1,...,z_n)=Res(z_1,...,z_n)\frac{g(z_1)...g(z_n)g(z_1^{-1}...z_n^{-1})}{g(z_1^k)...g(z_n^k)g(z_1^{-k}...z_n^{-k})~g(z_1^{k+1})...g(z_n^{k+1})g(z_1^{-k-1}...z_n^{-k-1})}.$$
It is clear that $H(z_1,...,z_n)$ is holomorphic, symmetric and it follows from (\ref{per}), (\ref{fun}) that $H(z_1,...,z_n)$ satisfies the identity
\begin{equation}\label{per1}
H(qz_1,...,z_n)=q^{k(k+1)}z_1^{2k^2+2k}z_2^{k^2+k}...z_n^{k^2+k}H(z_1,...,z_n). 
\end{equation}
However, one can check that any nonzero Laurent series in $z_1,...,z_n$ with these properties diverges. Therefore, $H$ (and also $Res$) is identically zero. $\Box$

{\bf Lemma 5.5.} The functions $f(z_1,...,z_k)=f_{k,k}(z_1,...,z_k)$ and $g(z)$ defined by (\ref{fun}) satisfy the equations (\ref{fah1}) and (\ref{hom}). Here $q\in\C$ and $|q|<1$.

{\bf Proof.} It follows from the previous Lemma that $f_{n,k}(z_1,...,z_k)$ can be written as convergent Laurent series for all $n,k=1,2,...$. However, any Laurent series satisfying (\ref{fpr}) with $n=k$ is homogeneous and satisfies (\ref{hom}). Furthermore, we have $f_{k+1,k}(z_1,...,z_{k+1})=0$ identically because any nonzero Laurent series satisfying (\ref{fpr}) with 
$n=k+1$ diverges. On the other hand, it follows from (\ref{req}) that the l.h.s. of (\ref{fah1}) is equal to $-f_{k+1,k}(z_1,...,z_{k+1})g(z_1^{-1}...z_{k+1}^{-1})=0.$ $\Box$

{\bf Theorem 5.2.} Let $I$ be an ideal of $A$ generated by $y_i,~i\in\Z$ where the generating function (\ref{fh1}) of the coefficients $a_{\alpha_1,...,\alpha_k}$ is equal to $f_{k,k}(z_1,...,z_k)$. Then $A/I$ is a flat deformation of the algebra $A/I_0$ where $I_0$ is the ideal generated by $x_ix_{i+1}...x_{i+k-1},~i\in\Z$.

{\bf Proof.} It follows from (\ref{fun}) that in the limit $q\to 0$ we have $g(z)=1-z$. Using recurrence relation (\ref{req}) one can prove by induction that in the limit $q\to 0$ we have
$$f_{n,k}(z_1,...,z_k)=\sum_{\substack{1\leq i_1,...,i_n\leq k,\\i_a\ne i_b \text{~if~}a\ne b}} z_1^{i_1}...z_n^{i_n}.$$
In particular, if $n=k$, then $\{i_1,...,i_k\}=\{1,2,...,k\}$ and we have in the limit $q\to0$
$$f(z_1,...,z_k)=f_{k,k}(z_1,...,z_k)=\sum_{\sigma\in S_k}x_{\sigma(1)}x_{\sigma(2)}^2...x_{\sigma(k)}^k=\psi_k(x_1...x_k).$$
More generally, it follows from (\ref{fih}) that in the limit $q\to0$ we have
$$f_i(z_1,...,z_k)=\psi_k(y_i)=\psi(x_ix_{i+1}...x_{i+k-1}).$$
Therefore, $A/I$ is a deformation of $A/I_0$. It follows from Theorem 4.2 that this deformation is flat iff relations (\ref{relh1}) hold for some sequence $b_{\beta}$. Moreover, it follows from Lemma 5.5 that these relations (\ref{relh1}) hold if generating functions for sequences $a_{\alpha_1,...,\alpha_k},b_{\beta}$ are given by (\ref{fun}). $\Box$

\subsection{Fermionic case}

Let $F_0=\Lambda^{*}[\xi_i,~i\in\Z]$ be generated by $\{\xi_i,~i\in\Z\}$ with defining relations $\xi_i\xi_j+\xi_j\xi_i=0,~i,j\in\Z$. It is a $\Z^2$-graded associative algebra, we assume 
$\text{deg}(\xi_i)=(i,1)$. Let $F$ be the projective limit of $F_0/(\xi_i,~|i|>N)$ in the category of $\Z^2$-graded associative algebras. Informally, $F$ is an extension of $F_0$, we allow certain infinite sums similarly as we did in Section 2. We have functional realization of $F$ in the space of formal Laurent series in several variables defined by
$$\psi_l(\xi_{i_1}...\xi_{i_l})=\sum_{\sigma\in S_l}\text{sign}(\sigma)z_{\sigma(1)}^{i_1}...z_{\sigma(l)}^{i_l}, $$
the only difference from the commutative case is that these Laurent series are skew-symmetric in the fermionic case. 

We are interested in flat deformations of quotients of the algebra $F$ in the class of skew-symmetric associative algebras. One can derive functional equations similar to (\ref{fa1}), (\ref{fah1}) as a 
criterion for flatness.

{\bf Conjecture 5.2.} Let $I_0$ be the ideal in the algebra $F$ generated by $\xi_i\xi_{i+1},~i\in\Z$. Let $I$ be the ideal in the algebra $F$ generated by
$$y_i=\xi_i\xi_{i+1}-q\xi_{i-1}\xi_{i+2}+q^3\xi_{i-2}\xi_{i+3}-...=\sum_{k=0}^{\infty}(-1)^k\xi_{i-k}\xi_{i+k+1}q^{\frac{k(k+1)}{2}},~~~i\in\Z.$$
Then the algebra $F/I$ is a flat deformation of the algebra $F/I_0$.

More generally, fix $k>1$. Let 
$$g(z)=(1-z)\prod_{i=1}^{\infty}(1-q^i)(1-q^iz)(1-q^iz^{-1})=\sum_{i\in\Z}(-1)^iz^iq^{\frac{i(i-1)}{2}}.$$

{\bf Conjecture 5.3.} Let $I_0$ be the ideal in the algebra $F$ generated by $\xi_i\xi_{i+1}...\xi_{i+k-1},~i\in\Z$. Let $I$ be the ideal in the algebra $F$ generated by $y_i~i\in\Z$ such that 
$$\psi_k(y_i)=z_1^iz_2^{i+1}z_3^{i+2}...z_k^{i+k-1}\prod_{1\leq i<j\leq k}g\Big(\frac{z_i}{z_j}\Big).$$
Then the algebra $F/I$ is a flat deformation of the algebra $F/I_0$.

{\bf Remark 5.4.} To support these conjectures we have verified the corresponding functional equations.

\section{Conclusion and outlook}

We have a strong feeling that algebras constructed in this paper constitute just a tip of an iceberg. In particular, commutative algebras studied in \cite{f1} should also admit this kind 
of ``elliptic'' deformations. Moreover, all these algebras can possibly be deformed further into associative non-commutative algebras. It would be interesting to understand if main objects 
of Representation Theory and Vertex Algebras (such as universal enveloping algebras of affine Lie algebras and Virasoro algebra) also admit similar elliptic deformations. 

The algebra $A$ admits cyclic reductions: we can set $x_{i+N}=x_i,~i\in\Z$ where $N$ is fixed. Elliptic deformations of quotients $A/I_0$ studied in this paper also admit such reductions. 
After cyclic reduction of these algebras we obtain families of commutative algebras with generators $x_i,~i\in\Z/N$. These algebras and their non-commutative deformations might be connected with Elliptic Algebras \cite{od}.

\bigskip

BF: National Research University, Higher School of Economics, Russia, Moscow, 101000, Myasnitskaya ul., 20

e-mail: bfeigin@gmail.com

\smallskip

AO: Brock University, Department of Mathematics and Statistics, 1812 Sir Isaac Brock Way, St. Catharines, ON, L2S 3A1 Canada

e-mail: aodesski@brocku.ca

\bigskip

This work has been funded by the  Russian Academic Excellence Project '5-100'.

\end{document}